\documentclass{article}
\catcode`\@=11
\newtheorem{theorem}{Theorem}
\newtheorem{lemma}{Lemma}
\newtheorem{proposition}{Proposition}
\newtheorem{definition}{Definition}
\newtheorem{example}{Example}
\newtheorem{remark}{Remark}

\newtheorem{corollary}{Corollary}
\def\demo{\noindent{\bf Proof .-}}
\def\section{\@startsection {section}{1}{\z@}{-3.5ex plus -1ex
minus-.2ex}{2.3ex plus .2ex}{\normalsize\bf}}

\begin{document}
\begin{center}
{\Large\bf \textsc{Arithmetical ranks of Stanley-Reisner ideals via linear algebra}}
\end{center}
\vskip.5truecm
\begin{center}
{Margherita Barile\footnote{Partially supported by the Italian Ministry of Education, University and Research.}\\ Dipartimento di Matematica, Universit\`{a} di Bari, Via E. Orabona 4,\\70125 Bari, Italy}\footnote{e-mail: barile@dm.uniba.it, Fax: 0039 080 596 3612}
\end{center}
\vskip1truecm
\noindent
{\bf Abstract} We present some examples of squarefree monomial ideals whose arithmetical rank can be computed using linear algebraic considerations.

\section*{Introduction and Preliminaries}
The {\it arithmetical rank} (ara) of an ideal $I$ in a commutative Noetherian ring $R$ is the minimal number $s$ of elements $a_1,\dots, a_s$ of $R$ such that $\sqrt I=\sqrt{(a_1,\dots, a_s)}$; one can express this equality by saying that $a_1,\dots, a_s$  {\it generate} $I$ {\it up to radical}. In general height\,$I\leq\,$ara\,$I$; if equality holds, $I$ is called a {\it set-theoretic complete intersection}. In this paper we determine the arithmetical ranks of some ideals which, in a polynomial ring over a field, are generated by monomials. Some results in this direction have already been proven in several works of the same author (\cite{B1}, \cite{B2}, \cite{B3}, \cite{B4}, \cite{B5}, \cite{B6}, \cite{B7}, \cite{B8}). In this note we study new examples in which the problem cannot be solved by means of the previously known methods; in particular we settle two cases which were left open in \cite{B8}.  The technique we are going to develop is based on linear algebraic considerations as in \cite{B1} and in \cite{B2}, but here, unlike in those papers, our approach is completely characteristic-free.\newline
Our results are best presented if the ideals are placed in  a combinatorial framework. First of all observe that, since the arithmetical rank of any ideal remains unchanged when the latter is replaced  by its radical, we can restrict our attention to ideals generated by squarefree monomials.   Let $X$ be a  finite set of indeterminates over the field $K$. A {\it simplicial complex} on $X$ is a set $\Delta$ of subsets of $X$ such that for all $x\in X$, $\{x\}\in\Delta$ and whenever $F\in\Delta$ and $G\subset F$, then $G\in\Delta$. The elements of $\Delta$ are called  {\it faces}, whereas $X$ is called the {\it vertex set} of $\Delta$, and the elements of $X$ are called the {\it vertices} of $\Delta$. If $\Delta$ consists of all subsets of its vertex set, then it is called a {\it simplex}. The simplicial complex $\Delta$ can be associated with an ideal $I_{\Delta}$ of the polynomial ring $R=K[X]$, which is generated by all monomials whose support is not a face of $\Delta$; $I_{\Delta}$ is called the {\it  Stanley-Reisner ideal} of $\Delta$ (over $K$). Its minimal monomial generators are the products of the elements of the minimal non-faces of $\Delta$, and these are squarefree monomials. In fact, this construction provides a one-to-one correspondence between the simplicial complexes on $X$ and the squarefree monomial ideals of $K[X]$.  We briefly recall some basic facts about Stanley-Reisner ideals, for which we refer to the extensive treatment given in \cite{BH}, Section 5.\par\smallskip\noindent
The minimal primes of $I_{\Delta}$ are the ideals of the form
$$P_F=(X\setminus F),\qquad\mbox{ where $F$ is any maximal face of $\Delta$}.$$
\noindent
It follows that the height of $I_{\Delta}$ is equal to $\vert X\vert-max_{F\in\Delta}\vert F\vert$. The number $d=max_{F\in\Delta}\vert F\vert-1$ is called the {\it dimension} of  $\Delta$; it is evidently equal to $\dim\, K[X]/I_{\Delta}-1$. The ideal $I_{\Delta}$ is unmixed if and only if all maximal faces of $\Delta$ have the same cardinality: we then say that $\Delta$ is {\it pure}. Note that the maximal faces of a pure one-dimensional  simplicial complex $\Delta$ form a graph; in general, the nonempty faces of a simplicial complex form a hypergraph, which we will use to represent $\Delta$ pictorially.   \par\smallskip\noindent 
If $I_{\Delta}$ is a set-theoretic complete intersection, then one can prove that it is Cohen-Macaulay, from which one concludes that it is unmixed, i.e.,  $\Delta$ is pure. If $I_{\Delta}$ is Cohen-Macaulay for any field $K$, then we will call $\Delta$ a {\it Cohen-Macaulay simplicial complex}. According to a well-known geometric characterization, a pure one-dimensional simplicial complex is Cohen-Macaulay if and only if the associated graph is connected.\par\smallskip\noindent
It is not known whether the Stanley-Reisner ideal of every Cohen-Macaulay simplicial complex is a set-theoretic complete intersection. The question is unsolved in general even for the special class of one-dimensional Gorenstein simplicial complexes, whose associated graphs are the cycle graphs $C_n$ for $n\geq 1$. The answer is (trivially) affirmative for $n=1,2,3,4$, since in all these cases $I_{\Delta}$ is a complete intersection. In \cite{B8} we established that   the answer is also affirmative for $n=5$, which we deduced from a general criterion based on the divisibility relations between the products of monomial generators. This criterion, indeed, enabled us to prove the set-theoretic complete intersection property for various one-dimensional pure simplicial complexes but, as we observed, it does not solve the problem for $n=6$. In this paper we will settle this case by a different approach: we will perform a direct computation which explicitly involves Cramer's Rule. Vanishing of determinants and proportionality conditions between rows of matrices are also applied to give a characteristic-free treatment of another example from \cite{B8}. 
\newline
Knowing that a certain Stanley-Reisner ideal is a set-theoretic complete intersection is not only interesting for itself: in our main result  we will show how this allows us to determine the arithmetical rank of various other related Stanley-Reisner ideals, whose simplicial complexes are derived from the original one through a simple geometric construction.\newline
The arithmetical rank of squarefree monomial ideals has also been intensively investigated by Lyubeznik (\cite{L1}, \cite{L2}) and Terai et al. (\cite{KTY}, \cite{T1}, \cite{T2}).\newline 
We recall the following result due to by Schmitt and Vogel, which will be useful in our proofs.
\begin{lemma}\label{Schmitt}{\rm (see \cite{SV}, p.~249)}. Let $P$ be a finite subset of elements of $R$. Let $P_0,\dots, P_r$ be subsets of $P$ such that
\begin{list}{}{}
\item[(i)] $\bigcup_{l=0}^rP_l=P$;
\item[(ii)] $P_0$ has exactly one element;
\item[(iii)] if $p$ and $p''$ are different elements of $P_l$ $(0<l\leq r)$ there is an integer $l'$ with $0\leq l'<l$ and an element $p'\in P_{l'}$ such that $pp''\in(p')$.
\end{list}
\noindent
Let $0\leq l\leq r$, and, for any $p\in P_l$, let $e(p)\geq1$ be an integer. We set $q_l=\sum_{p\in P_l}p^{e(p)}$. We will write $(P)$ for the ideal of $R$ generated by the elements of $P$.  Then we get
$$\sqrt{(P)}=\sqrt{(q_0,\dots,q_r)}.$$
\end{lemma} 
\section{The main theorem}
Without loss of generality, we shall throughout assume that the field $K$ is algebraically closed. \newline
We will consider the following two sets of vertices/indeterminates over $K$:  $X=\{x_1,\dots, x_n\}$ and $X'=X\cup\{x_0\}$. We set $R=K[X]$ and $R'=K[X']$. 
\begin{definition}\label{def1}{\rm  Let $F$ be any nonempty subset of $X$. We will call {\it cone from $x_0$ over F}, denoted co\,$_{x_0}F$, the simplex on the vertex set $F\cup\{x_0\}$. 
}
\end{definition}
\par\smallskip\noindent
Let $\Delta$ be a simplicial complex on the vertex set $X$, and let $F$ be any maximal face of $\Delta$. Then $\Delta'=\Delta\cup{\rm co}\,_{x_0}F$ is a simplicial complex on the vertex set $X'$.  
As an immediate consequence of Definition \ref{def1} we have the following
\begin{proposition}\label{primes} 
The maximal faces of $\Delta'$ are the maximal faces $G\ne F$ of $\Delta$ and $F\cup\{x_0\}$. Correspondingly, the minimal primes of $I_{\Delta'}$ are the ideals $P_G+(x_0)$ and $P_F$. 
\end{proposition}
\par\smallskip\noindent
We are now ready to state our main result, which shows a relation between the arithmetical ranks of the ideal $I_{\Delta}$ of $R$ and of the ideal $I_{\Delta'}$ of $R'$.
\begin{theorem}\label{main} Suppose that $I_{\Delta}$ is a set-theoretic complete intersection. Then 
$${\rm ara}\,I_{\Delta'}=\,{\rm ara}\,I_{\Delta}+1.$$ 
\end{theorem}
\demo  Let $t=\,$height\,$I_{\Delta}$. Since $I_{\Delta}$ is unmixed,  we have that $P_F=(x_{i_1},\dots, x_{i_t})$ for some (pairwise distinct) indices $i_1,\dots, i_t$.  We may assume that $i_j=j$ for all $j=1,\dots, t$. By assumption  $t=\,$ara\,$I_{\Delta}$, i.e., there are $q_1,\dots, q_t\in R$ which generate $I_{\Delta}$ up to radical. We will use these elements to construct $t+1$ elements of $R'$ which generate $I_{\Delta'}$ up to radical. This will allow us to conclude that ara\,$I_{\Delta'}\leq t+1$. The opposite inequality will be proved right afterwards.\newline
  We set $J=(q_1,\dots, q_t)$. For all indices $i=1,\dots, t$, we have that $q_i\in I_{\Delta}$, whence $q_i\in (x_1,\dots, x_t)$, i.e., $q_i=\sum_{j=1}^t a_{ij}x_j$ for some $a_{ij}\in R$. Recall that a polynomial belongs to a given monomial ideal   if and only if each of its monomial terms is  divisible by some monomial generator of this ideal. Therefore, up to eliminating redundant terms, we may assume that 
\begin{equation}\label{0}\qquad\qquad a_{ij}x_j\in I_{\Delta}\qquad\qquad\qquad\mbox{ for all indices } i,j\in\{1,\dots, t\}.\end{equation}
\noindent
Consider the following ring homomorphism:
$$\phi:R\longrightarrow R$$
$$f(x_1,\dots, x_n)\mapsto f(x_1^2,\dots, x_n^2),$$
\noindent
and let $\bar J =\phi(J)$. We now show that
\begin{equation}\label{1} \sqrt{\bar J}=I_{\Delta}.\end{equation}
\noindent
For all $i=1,\dots, t$, set $\bar q_i=\phi(q_i)$, so that 
\begin{equation}\label{J}\bar J=(\bar q_1,\dots, \bar q_t).\end{equation}
\noindent
First of all note that, since all monomial terms of $q_i$ belong to $I_{\Delta}$, the same is true for all monomial terms of $\bar q_i$, because these are the squares of the monomial terms of $q_i$. Hence $\bar J\subset I_{\Delta}$, which implies that $\sqrt{\bar J}\subset I_{\Delta}$. Next we prove the opposite inclusion. Let $g$ be any monomial generator of $I_{\Delta}$. Since, by assumption, $\sqrt{J}=I_{\Delta}$, for some positive integer $a$ we have that $g^a\in J$, i.e., there are $f_1,\dots, f_t\in R$ such that $g^a=\sum_{i=1}^tf_iq_i$. But then 
$$g^{2a}=\phi(g^a)=\sum_{i=1}^t\phi(f_i)\bar q_i,$$
\noindent
which shows that $g\in\sqrt{\bar J}$. Thus $I_{\Delta}\subset\sqrt{\bar J}$, which completes the proof of (\ref{1}).\newline For all $i,j\in\{1,\dots, t\}$, set $\bar a_{ij}=\phi(a_{ij})x_j$, so that, for all $i=1,\dots, t$, 
\begin{equation}\label{barq}\bar q_i=\sum_{j=1}^t\bar a_{ij}x_j.\end{equation}
\noindent
The monomial terms of every $\bar a_{ij}$ are all of the form $b^2x_j$, where $b$ is a monomial term of $a_{ij}$; hence, by (\ref{0}), $bx_j\in I_{\Delta}$, so that $b^2x_j\in I_{\Delta}$, and, for all $i,j\in\{1,\dots, t\}$,  
\begin{equation}\label{2}\bar a_{ij}\in I_{\Delta},\qquad\mbox{whence}\qquad \bar q_i\in I_{\Delta}.\end{equation}
\noindent
 Let $\bar A=(\bar a_{ij})_{i,j=1,\dots, t}$, and set $A'=\bar A + x_0Id_t$, where $Id_t$ denotes the $t\times t$ identity matrix. Moreover, let 
\begin{equation}\label{det} D=\det A'-x_0^t.\end{equation}
\noindent
 Now, by definition of determinant,  $D$ is the sum of products each of which involves at least one entry of $\bar A$ as a factor; in view of (\ref{2}), it follows that  
\begin{equation}\label{D} D\in I_{\Delta}.\end{equation}
\noindent
 Set 
\begin{equation}\label{J'}J'=(D,\bar q_1+x_0x_1,\dots, \bar q_t+x_0x_t).\end{equation}
\noindent
We claim that
\begin{equation}\label{3} \sqrt{J'}=I_{\Delta'}.\end{equation}
\noindent
Note that, by definition of Stanley-Reisner ideal,  
\begin{equation}\label{4} I_{\Delta'}=I_{\Delta}R'+(x_0x_1,\dots, x_0x_t).\end{equation}
\noindent
Now (\ref{2}), (\ref{D}), (\ref{J'}) and  (\ref{4}) imply that  $J'\subset I_{\Delta'}$, so that  $\sqrt{J'}\subset I_{\Delta'}$. For the opposite inclusion we use Hilbert's Nullstellensatz. Let ${\bf x}\in K^{n+1}$ be such that all elements of $J'$ vanish at ${\bf x}$. We show that ${\bf x}$ annihilates all elements of $I_{\Delta'}$. In the rest of the proof, we shall identify each polynomial with its value at ${\bf x}$. Thus our assumption can be formulated in the form:
\begin{eqnarray}\label{equationa} D&=&0,\\
\bar q_1+x_0x_1&=&\cdots=\bar q_t+x_0x_t=0.\label{equationb}
\end{eqnarray}
\noindent
 We distinguish between two cases. First suppose that $\det A'\ne 0$. Note that $A'$ is the matrix of coefficients of the square system of homogeneous linear equations   
$$\sum_{i=1}^t(\bar a_{ij}+\delta_{ij}x_0)y_j=0\qquad(j=1,\dots, t)$$
in the unknowns $y_1,\dots, y_t$. By Cramer's Rule it only has  the trivial solution. Therefore, in view of (\ref{barq}), (\ref{equationb}) implies that $x_1=\cdots=x_t=0$. But, in view of (\ref{4}), $I_{\Delta'}\subset (x_1,\dots, x_t)$, so that ${\bf x}$ annihilates all elements of $I_{\Delta'}$. Now suppose that $\det A'=0$. Then, in view of (\ref{det}), from (\ref{equationa}) we deduce that $x_0=0$, so that from (\ref{equationb}) we further have that $\bar q_1=\cdots=\bar q_t=0$. These equalities, together with (\ref{1}) and (\ref{J}), imply that all elements of $I_{\Delta}$ vanish at ${\bf x}$.  In view of (\ref{4}), we again conclude that ${\bf x}$ annihilates all elements of $I_{\Delta'}$. This completes the proof of (\ref{3}). From (\ref{3}) we deduce that ara\,$I_{\Delta'}\leq t+1$, as required. On the other hand, we have that ara\,$I_{\Delta'}$ is greater than or equal to the height of all minimal primes of $I_{\Delta'}$: this is a consequence of Krull's Principal Ideal Theorem (see \cite{M}, Theorem 13.5). In view of Proposition \ref{primes} we thus conclude that 
ara\,$I_{\Delta'}\geq t+1$. Hence ara\,$I_{\Delta'}= t+1$, as was to be shown, and $D,\bar q_1+x_0x_1,\dots, \bar q_t+x_0x_t$ are  $t+1$ elements of $R'$ generating   $I_{\Delta'}$ up to radical. This completes the proof.
\par\smallskip\noindent
\begin{example}\label{example1}{\rm
Consider the simplicial complex $\Delta$ on the vertex set $\{x_1,\dots, x_5\}$ whose maximal faces are $F=\{x_1, x_2\}$, $\{x_2, x_3\}$, $\{x_3, x_4\}$, $\{x_4, x_5\}$, $\{x_5, x_1\}$. It is associated with the cycle graph $C_5$.
Its (Gorenstein) Stanley-Reisner ideal in the polynomial ring $R=K[x_1,\dots, x_5]$ is 
$$I_{\Delta}=(x_1x_3,\ x_1x_4,\ x_2x_4,\ x_2x_5,\ x_3x_5), $$
\noindent
whose minimal prime decomposition  is
$$I=(x_3, x_4, x_5)\cap(x_1, x_4, x_5)\cap(x_1, x_2, x_3)\cap(x_2, x_3, x_4)\cap(x_1, x_2, x_5).$$
\noindent
In \cite{B8}, Example 1, we proved that it is a set-theoretic complete intersection, by showing that it is generated, up to radical, by the following three elements:
\begin{equation}\label{q} q_1=x_1x_3,\ q_2=x_1x_4+x_2x_5,\ q_3=x_2x_4+x_3x_5.\end{equation} 
\noindent
According to Proposition \ref{primes}, the corresponding simplicial complex $\Delta'=\Delta\cup\,{\rm co}\,_{x_0}F$ on the vertex set $\{x_0, x_1,\dots, x_5\}$ has the following maximal faces: ${\rm co}\,_{x_0}F=\{x_1, x_2, x_0\}$, $\{x_2, x_3\}$, $\{x_3, x_4\}$, $\{x_4, x_5\}$, $\{x_5, x_1\}$.
The monomial generators of $I_{\Delta'}$ are:
$$x_1x_3, x_1x_4, x_2x_4, x_2x_5, x_3x_5, x_0x_3, x_0x_4, x_0x_5.$$
\noindent
According to Theorem \ref{main}, ara\,$I_{\Delta'}=4$; we construct four elements generating $I_{\Delta'}=4$ up to radical applying the procedure described in the proof of the theorem to the elements $q_1, q_2, q_3$ presented in (\ref{q}). 
With respect to the notation introduced above, we have that $P_F=(x_3, x_4, x_5)$ (i.e., in the proof of Theorem \ref{main}, $i_1=3$, $i_2=4$, $i_3=5$),  and the matrix $(a_{ij})_{i,j=1,\dots,3}$ is
$$A=\left(\begin{array}{ccc}
x_1&0&0\\
0&x_1&x_2\\
0&x_2&x_3
\end{array}\right),
$$
\noindent
from which we derive the matrix
 $$\bar A=\left(\begin{array}{ccc}
x_1^2x_3&0&0\\
0&x_1^2x_4&x_2^2x_5\\
0&x_2^2x_4&x_3^2x_5
\end{array}\right)
$$
\noindent
and finally, the matrix
$$
A'=\bar A +x_0Id_3=
\left(\begin{array}{ccc}
x_1^2x_3+x_0&0&0\\
0&x_1^2x_4+x_0&x_2^2x_5\\
0&x_2^2x_4&x_3^2x_5+x_0
\end{array}\right).
$$
\noindent
It follows that $I_{\Delta'}$ is generated up to radical by the following four elements:
\begin{eqnarray*}&&D=\det A'-x_0^3=x_1^4x_3^3x_4x_5+x_0x_1^2x_3^2x_4x_5+x_0x_1^2x_3^3x_5+x_0^2x_3^2x_5\\
&&\qquad\qquad\qquad\quad+x_0x_1^4x_3x_4+x_0^2x_1^2x_4 +x_0^2x_1^2x_3-x_1^2x_2^4x_3x_4x_5-x_0x_2^4x_4x_5,\\
&&x_1^2x_3^2+x_0x_3,\\
&&x_1^2x_4^2+x_2^2x_5^2+x_0x_4,\\
&&x_2^2x_4^2+x_3^2x_5^2+x_0x_5.
\end{eqnarray*}
\noindent}
\end{example}
\pagebreak
\section{More examples of computations of arithmetical ranks via linear algebra}
In this section we present a  class of squarefree monomial ideals whose arithmetical ranks can be determined using determinants and Cramer's Rule.\par\smallskip\noindent
For all integers $n\geq6$ let $I_n$ be the ideal of $R=K[x_1,\dots, x_n]$ generated by the following squarefree monomials:
\begin{eqnarray}\label{generators} &&x_1x_3,\dots, x_1x_{n-1},\nonumber\\
&&x_2x_4,\dots, x_2x_n,\nonumber\\
&&x_3x_5,\dots, x_3x_n,\nonumber\\
&&x_4x_n,\dots, x_{n-2}x_n.
\end{eqnarray}
\noindent
Our aim is to determine the arithmetical rank of $I_n$ and to give ara\,$I_n$ elements of $R$ generating $I_n$ up to radical. 
To this end we introduce the matrix $B=(b_{ij})_{i,j=1,\dots, n-3}$ whose entries are defined as follows.
\begin{list}{}{}
\item{I.}
\begin{list}{}{}
\item{(i)} $b_{jj}=x_1$ for $j=1,\dots, n-4$;
\item{(ii)} $b_{n-3, n-3}=x_3$.
\end{list}
(All elements on the main diagonal are equal to $x_1$, except the last one, which is $x_3$.)
\item{II.} $b_{j+1, j}=x_2$ for $j=1,\dots, n-4$.\newline
(All elements on the first ``under-diagonal'' are equal to $x_2$.)
\item{III.} $b_{ij}=0$ if $i\geq j+2$.\newline
(All elements below the first ``under-diagonal'' are equal to 0.)
\item{IV.}  
 $b_{j-1, j}=x_3x_{3+j}$ for $j=2,\dots, n-4$.\newline
(The elements on the first ``over-diagonal'', except the last one, are  equal to $x_3x_5, \dots, x_3x_{n-1}$, respectively.)
\item{V.} 
\begin{list}{}{}
\item{(i)} $b_{1, n-3}=x_2$;
\item{(ii)} $b_{i, n-3}=x_{2+i}$ for $i=2,\dots, n-4$.\newline
\end{list}
(The first element of the last column is $x_2$, the following ones, except the last one, are $x_4$, \dots, $x_{n-2}$, respectively.)
\item{VI.} $b_{ij}=0$ if $j\geq i+2$ and $j\neq n-3$.\newline
(All elements above the first ``over-diagonal'', except those in the last column, are equal to 0.)
\end{list}
\noindent
Hence 
$$B=\left(\begin{array}{ccccccc}
x_1&x_3x_5&&&&&x_2\\
x_2&x_1&x_3x_6&&&&x_4\\
&x_2&x_1&\ddots&&&x_5\\
&&x_2&\ddots&\ddots&&\vdots\\
&&&\ddots&\ddots&x_3x_{n-1}&\vdots\\
&&&&\ddots&x_1&x_{n-2}\\
&&&&&x_2&x_3
\end{array}
\right),
$$
\noindent
where it is understood that the empty triangular regions are occupied by zeros. \par\smallskip\noindent
Set 
\begin{equation}\label{B}D=\det B-(-1)^nx_2^{n-3}.\end{equation}
\begin{lemma}\label{lemma}The following hold:
\begin{list}{}{}
\item{(a)} $D\in I_n$;
\item{(b)} $D-x_1^{n-4}x_3\in x_2(x_4,\dots, x_n)$.
\end{list} 
\end{lemma}
\demo
Let $p$ be a nonzero term in the Laplace expansion of $\det B$ with respect to the first columun of $B$. According to the above prescriptions I.(i), II., and III., the only nonzero entries of the first column of $B$ are $b_{11}=x_1$ and $b_{21}=x_2$.  Hence $p$ is the product of factors obtained by picking $x_1$ or $x_2$ in the first column and one entry in each of the columns of the submatrix of $B$ obtained by omitting the first row and the first column, or the second row and the first column, respectively. \newline
For the proof of (a), we show that one of the terms $p$ is $(-1)^nx_2^{n-3}$, whereas the others are all divisible by some of the generators of $I_n$ listed in (\ref{generators}). First suppose that $x_1$ is the entry picked in the first column. Then the entry picked in the $(n-3)$-th column is one of $b_{2,n-3}=x_4,\dots, b_{n-4, n-3}=x_{n-2}$, $b_{n-3, n-3}=x_3$.  Hence $p$ is divisible by one of the following generators of $I_n$: $x_1x_4,\dots, x_1x_{n-2}, x_1x_3$, as required. Now suppose that $x_2$ is the entry picked in the first column. Also assume that $b_{j+1,j}=x_2$ is the entry picked in all columns with the indices $j=1,\dots, k$, for some $k$ such that $1\leq k\leq n-4$.  Suppose that $k$ is maximal with respect to these conditions. If $k=n-4$, then the entry picked in the $(n-3)$-th column is $b_{1,n-3}=x_2$, and then $p=(-1)^nx_2^{n-3}$. So suppose that $k<n-4$. Then the entry picked in the $(k+1)$-th column is $b_{i,k+1}$, with $k+1\ne n-3$, for some $i\not\in\{2,\dots,k+1\}$; by maximality of $k$, it cannot be $b_{k+2, k+1}=x_2$. Moreover, for all $i\geq k+3$, by III.  we have that $b_{i,k+1}=0$. Thus   the entry picked in the $(k+1)$-th column is $b_{1, k+1}$, which, in view of I.(i), IV. and VI., is nonzero only for $k=1$. Hence this entry is $b_{12}=x_3x_5$. Therefore, $p$  is divisible by $x_3x_5$, which is one of the generators of $I_n$. This completes the proof of (a).
\newline
We now prove (b). We first show that if $x_2$ divides $p\ne (-1)^nx_2^{n-3}$, then one of $x_4,\dots, x_n$ divides $p$. We have seen in the last part of the proof of (a) that $x_5$ divides $p$ whenever $p$ is obtained by picking $b_{21}=x_2$ in the first column. So assume that $p$ is divisible by $x_2$ and $b_{11}=x_1$ is the entry picked in the first column. Then $x_2=b_{j+1,j}$ is the entry picked in the $j$-th column for some $j\in\{2,\dots, n-4\}$. Let $j$ be minimal with respect to this property. Then the entry picked in the $j$-th row is not $b_{j,j-1}=x_2$ (by minimality), nor $b_{jj}=x_1$, which lies in the $j$-th column. Hence, discarding zero entries, this entry is necessarily $b_{j, j+1}=x_3x_{4+j}$ if $j\ne n-4$ or $b_{n-4, n-3}=x_{n-2}$ if $j=n-4$.  But then $x_{3+j}$ or $x_{n-2}$ divides $p$, as required. 
There remains to show that, if $p$ is not divisible by $x_2$, then $p=x_1^{n-4}x_3$. So assume that $p$ is not divisible by $x_2$. Then it has been obtained by picking $x_1$ in the first column. Suppose that $b_{jj}=x_1$ is the entry picked in all the columns with the indices $j=1,\dots, k$, for some $k$ such that $1\leq k\leq n-4$. Then the entry picked in the $(k+1)$-th column is $b_{i,k+1}$ for some $i\not\in\{1,\dots,k\}$, and it cannot be $b_{k+2, k+1}=x_2$. According to III., the only nonzero available entry is then $b_{k+1, k+1}$. This shows by finite induction that $p=\prod_{j=1}^{n-3}b_{jj}=x_1^{n-4}x_3$, as required. This completes the proof of (b), and of the Lemma. 
\par\bigskip\noindent
Let
\begin{equation}\label{star}q_i=\sum_{j=1}^{n-3}b_{ij}x_{3+j}\qquad\mbox{ for }i=1,\dots, n-3.\end{equation}
\noindent
We can now prove the following result. 
\begin{proposition}\label{proposition} For all integers $n\geq 6$, 
$$I_n=\sqrt{(D, q_1,\dots, q_{n-3})}.$$
\end{proposition}
\demo Set $J=(D, q_1,\dots, q_{n-3})$. For all $i,j=1\,\dots, n-3$ set $s_{ij}=b_{ij}x_{3+j}$. We first show that the $s_{ij}$'s (i.e., the summands of the $q_i$'s) are, up to repeated factors, the generators of $I_n$ listed in (\ref{generators}), with the only exception of $x_1x_3$. In fact these terms $s_{ij}$ are:
\begin{eqnarray*}
s_{14}&=&b_{11}x_4=x_1x_4,\dots,s_{n-4,n-4}=b_{n-4,n-4}x_{n-1}=x_1x_{n-1},\quad\mbox{(see I.(i))}\\
 s_{21}&=&b_{21}x_4=x_2x_4,\dots,s_{n-3,n-4}=b_{n-3,n-4}x_{n-1}=x_2x_{n-1},\quad\mbox{(see II.)}\\
&&\qquad\qquad\qquad\quad \qquad s_{1,n-3}=b_{1,n-3}x_n=x_2x_n,\qquad\qquad\mbox{(see V.(i))}\\
s_{12}&=&b_{12}x_5=x_3x_5^2\dots,s_{n-5,n-4}=b_{n-5,n-4}x_{n-1}=x_3x_{n-1}^2,\quad\mbox{(see IV.)}\\
&&\qquad\qquad\qquad\quad \qquad s_{n-3,n-3}=b_{n-3,n-3}x_n=x_3x_n,\quad\ \ \mbox{(see I.(ii))}\\
s_{2,n-3}&=&b_{2,n-3}x_n=x_4x_n,\dots,s_{n-4,n-3}=b_{n-4,n-3}x_n=x_{n-2}x_n.\ \mbox{(see V.(ii))}
\end{eqnarray*}
\noindent
It follows that $q_1,\dots, q_{n-3}\in I_n$. In view of Lemma \ref{lemma}(a) this implies that $J\subset I_n$, whence $\sqrt J\subset I_n$. 
\newline
For the proof of the opposite inclusion we use Hilbert's Nullstellensatz. Let ${\bf x}\in K^n$ be such that $D, q_1,\dots, q_{n-3}$ vanish at ${\bf x}$. We show that ${\bf x}$ annihilates all generators of $I_n$ listed in (\ref{generators}). In the rest of the proof, unless otherwise indicated, we shall identify every polynomial with its value at ${\bf x}$. Hence our assumption can be stated as follows:
\begin{eqnarray}
D&=&0\label{D0},\\
q_1&=&\cdots=q_{n-3}=0.\label{qi}
\end{eqnarray}
\noindent 
We distinguish between two cases. First assume that $\det B\ne 0$. Note that $B$ is the matrix of the coefficients of the square system of homogeneous linear equations
$$\sum_{j=1}^{n-3}b_{ij}y_j=0\qquad\qquad(i=1,\dots, n-3)$$
in the unknowns $y_1,\dots, y_{n-3}$. By Cramer's Rule it follows that it  only has the trivial solution. Hence, in view of (\ref{star}), (\ref{qi}) implies that $x_4=\cdots =x_n=0$. This, together with Lemma \ref{lemma}(b), gives $D-x_1^{n-4}x_3=0$. From (\ref{D0}) it then follows that $x_1x_3=0$. Thus all generators of $I_n$ listed in (\ref{generators}) vanish at ${\bf x}$. Now assume that $\det B =0$. By (\ref{B}) and (\ref{D0}) we then have $D=x_2=0$, which, in view of Lemma \ref{lemma}(b), yields $x_1x_3=0$. There remains to prove that the monomials 
\begin{equation}\label{monomials} x_1x_4,\dots, x_1x_{n-1}, x_3x_5,\dots, x_3x_n, x_4x_n, \dots, x_{n-2}x_n
\end{equation}
\noindent
 vanish at ${\bf x}$ as well. For all $i=1,\dots, n-3$, let $\tilde q_i$ be the polynomial of $R$ obtained by setting $x_2=0$ in $q_i$. From the first part of the proof we know that the nonzero summands $s_{ij}$ of the $\tilde q_i$'s are, up to repeated factors, the monomials listed in (\ref{monomials}). Hence it suffices to show that Lemma \ref{Schmitt} can be applied to the following sequence of polynomials
$$p=x_1x_3, p_0=\tilde q_{n-3}=x_3x_n, p_1=\tilde q_1, \dots, p_{n-4}=\tilde q_{n-4}.$$
\noindent
We show that, for all $i=1,\dots, n-4$, the product of any two nonzero terms  $s_{ij}, s_{ik}$ of $\tilde q_i$ is divisible by $x_1x_3$ or $x_3x_n$ or some term $s_{i'j}$ with $1\leq i'<i$. For $i=1$ the only product to be considered is
$$s_{11}s_{12}=b_{11}x_4b_{12}x_5=x_1x_4x_3x_5^2,$$
\noindent
which is divisible by $x_1x_3$. For $i\in\{2,\dots, n-4\}$ we have to consider three different kinds of products (in the first and last we assume that $i\ne n-4$):
\begin{eqnarray*}
s_{ii}s_{i, i+1}&=&b_{ii}x_{3+i}b_{i,i+1}x_{4+i}=x_1x_{3+i}x_3x_{4+i}^2,\mbox{ which is divisible by }x_1x_3,\\
s_{ii}s_{i, n-3}&=&b_{ii}x_{3+i}b_{i,n-3}x_n=x_1x_{3+i}x_{2+i}x_n,\\
&&\qquad\mbox{ which is divisible by }s_{i-1,i-1}=b_{i-1,i-1}x_{3+i-1}=x_1x_{2+i},\\
s_{i,i+1}s_{i, n-3}&=&b_{i,i+1}x_{4+i}b_{i,n-3}x_n=x_3x_{4+i}^2x_{2+i}x_n,\mbox{ which is divisible by }x_3x_n.
\end{eqnarray*}
\noindent
This completes the proof.
\begin{corollary}\label{corollary1}
For all integers $n\geq 6$,
$${\rm ara}\,I_n=n-2.$$
\end{corollary}
\demo
From Proposition \ref{proposition} we have that ara\,$I\leq n-2$. On the other hand, in view of (\ref{generators}), $P=(x_1,\dots, x_{n-2})$ is a minimal prime of $I_n$. Hence by \cite{M}, Theorem 13.5, $n-2\leq\,$ara\,$I_n$. This completes the proof. \newline
\begin{example}\label{example2}{\rm The ideal $I_6$ of $K[x_1,\dots, x_6]$ is generated by the following squarefree monomials:
\begin{eqnarray*} &&x_1x_3, x_1x_4, x_1x_5,\\
&&x_2x_4, x_2x_5, x_2x_6, \\
&&x_3x_5, x_3x_6, \\
&&x_4x_6.
\end{eqnarray*}
\noindent
According to Proposition \ref{proposition} it is generated up to radical by the following four elements:
\begin{eqnarray*}
D&=&x_1^2x_3- x_1x_2x_4 -x_2x_3^2x_5,\\
q_1&=&x_1x_4+x_3x_5^2+x_2x_6,\\
q_2&=&x_2x_4+x_1x_5+x_4x_6,\\
q_3&=&x_2x_5+x_3x_6.
\end{eqnarray*}
\noindent
Note that $I_6$ is the Stanley-Reisner ideal of the  (Gorenstein) simplicial complex on the vertex set $\{x_1,\dots, x_6\}$ whose maximal faces are $\{x_1, x_2\}$, $\{x_2, x_3\}$, $\{x_3, x_4\}$, $\{x_4, x_5\}$, $\{x_5, x_6\}$, $\{x_6, x_1\}$. It is associated with the cycle graph $C_6$. \newline
The question, raised by Marcel Morales, whether the Stanley-Reisner rings of the Gorenstein simplicial complexes associated with the cycle graphs $C_n$  are all set-theoretic complete intersections remains open; in fact our method cannot allow us to settle the problem for $n\geq 7$.}  
\end{example}
\begin{example}\label{example3}{\rm Consider the ideal $I_7$ of $K[x_1,\dots, x_7]$, which  is generated by the following squarefree monomials:
\begin{eqnarray*} &&x_1x_3, x_1x_4, x_1x_5, x_1x_6,\\
&&x_2x_4, x_2x_5, x_2x_6, x_2x_7,\\
&&x_3x_5, x_3x_6, x_3x_7,\\
&&x_4x_7, x_5x_7.
\end{eqnarray*}
\noindent
According to Corollary \ref{corollary1}, ara\,$I_7=5$ and, by Proposition \ref{proposition}, $I_7$ it is generated up to radical by the following five elements:
\begin{eqnarray*}
D&=&x_1^3x_3-x_1^2x_2x_5-x_1x_2x_3^2x_6+x_1x_2^2x_4-x_1x_2x_3^2x_5+x_2^2x_3x_5^2,\\
q_1&=&x_1x_4+x_3x_5^2+x_2x_7,\\
q_2&=&x_2x_4+x_1x_5+x_3x_6^2+x_4x_7,\\
q_3&=&x_2x_5+x_1x_6+x_5x_7,\\
q_4&=&x_2x_6+x_3x_7.
\end{eqnarray*}
\noindent
Note that $I_7$ is the Stanley-Reisner ideal of the simplicial complex on the vertex set $\{x_1,\dots, x_7\}$ whose maximal faces are $\{x_1, x_2\}$, $\{x_2, x_3\}$, $\{x_3, x_4\}$, $\{x_4, x_5, x_6\}$, $\{x_6, x_7\}$, $\{x_7, x_1\}$.
}
\end{example}
We now complete the statement contained in Corollary \ref{corollary1}.
\begin{corollary}\label{corollary2} 
The ideal $I_n$ is a set-theoretic complete intersection if and only if $n=6$.
\end{corollary}
\demo Since, in view of (\ref{generators}), both $P=(x_1,\dots, x_{n-2})$ and $Q=(x_1, x_2, x_3, x_n)$ are minimal primes of $I_n$, $I_n$ is not unmixed if $n>6$, hence it is not a set-theoretic complete intersection. In view of Example \ref{example2}, this completes the proof.
\par\bigskip\noindent
\begin{remark}{\rm  Note that, up to renaming vertices, the simplicial complex in Example \ref{example3} can be obtained from the one in Example \ref{example2} by the cone construction described in Section 1. In particular, five elements generating $I_7$ up to radical can also be obtained from the four elements generating $I_6$ up to radical by means of Theorem \ref{main}.}   
\end{remark}
\section{One more set-theoretic complete intersection via linear algebra}
Next we present an example of Stanley-Reisner ideal whose set-theoretic complete intersection property can be shown by arguments which, in addition to Cramer's Rule, involve considerations on the proportionality between rows of a matrix. 
\begin{example}\label{example4}{\rm
In the polynomial ring $R=K[x_1,\dots, x_6]$ consider the ideal $I$ generated by the following squarefree monomials:
\begin{equation}\label{generators2}x_1x_4, x_1x_5, x_1x_2x_3, x_2x_4, x_2x_5, x_2x_6, x_3x_5, x_3x_6, x_4x_6.\end{equation}
\noindent
It is the Stanley-Reisner ideal of the Cohen-Macaulay simplicial complex on the vertex set $\{x_1,\dots, x_6\}$ whose maximal faces are $\{ x_1, x_2\}$, $\{ x_2, x_3\}$, $\{ x_1, x_3\}$, $\{ x_3, x_4\}$, $\{ x_4, x_5\}$, $\{ x_5, x_6\}$, $\{ x_6, x_1\}$.  The ideal $I$ is unmixed of height 4. 
Set
$$C=\left(\begin{array}{ccc}
x_1&x_2&x_3\\
x_2&x_3&x_4\\
0&x_1&x_2
\end{array}
\right),
$$
\noindent
so that
\begin{equation}\label{detC}\det C=2x_1x_2x_3-x_1^2x_4-x_2^3.\end{equation}
\noindent
We show that $I$ is generated, up to radical, by the following four elements:
\begin{eqnarray*}
D=\det C-x_1x_2x_3+x_2^3&=&x_1x_2x_3-x_1^2x_4,\\
q_1&=&x_1x_4+x_2x_5+x_3x_6,\\
q_2&=&x_2x_4+x_3x_5+x_4x_6,\\
q_3&=&x_1x_5+x_2x_6.
\end{eqnarray*}
\noindent
Since these all belong to $I$, it suffices to show that whenever ${\bf x}\in K^6$ annihilates them all, ${\bf x}$ annihilates all monomial generators of $I$, which are listed in (\ref{generators2}), as well. As usual, from now on we will identify all polynomials with their value at ${\bf x}$. So assume that
\begin{eqnarray}
x_1x_2x_3-x_1^2x_4&=&0,\label{equation1}\\
x_1x_4+x_2x_5+x_3x_6&=&0,\label{equation2}\\
x_2x_4+x_3x_5+x_4x_6&=&0,\label{equation3}\\
x_1x_5+x_2x_6&=&0.\label{equation4}
\end{eqnarray}
\noindent
Note that $C$ is the matrix of coefficients of the square system of homogeneous linear equations
\begin{eqnarray*}
x_1y_1+x_2y_2+x_3y_3&=&0,\\
x_2y_1+x_3y_2+x_4y_3&=&0,\\
x_1y_2+x_2y_3&=&0.
\end{eqnarray*}
 in the unknowns $y_1, y_2, y_3$. We distinguish between several cases and subcases.
\par\smallskip\noindent
\underline{Case 1}: Suppose that $\det C\ne 0$. Then, by Cramer's Rule, the above system of linear equations only has the trivial solution. Thus equalities (\ref{equation2}), (\ref{equation3}) and (\ref{equation4}), imply that $x_4=x_5=x_6$. From (\ref{equation1}) it then follows that $x_1x_2x_3=0$. Hence all generators listed in (\ref{generators2}) vanish at ${\bf x}$ in this case.
\par\smallskip\noindent
\underline{Case 2}: Suppose that $\det C=0$. Then (\ref{detC}) and (\ref{equation1}) imply
$$0=x_1x_2x_3-x_2^3=x_2(x_1x_3-x_2^2).$$
\noindent
Hence $x_2=0$ or $x_1x_3-x_2^2=0$.
\par\smallskip\noindent
\underline{Case 2.1}: Suppose that $x_2=0$, so that $x_1x_2x_3=x_2x_4=x_2x_5=x_2x_6=0$. Hence, from (\ref{equation4}) we get $x_1x_5=0$. Moreover, from (\ref{equation1}) we obtain $x_1x_4=0$, so that, by (\ref{equation2}), we also have that $x_3x_6=0$.  Now the only surviving summands in (\ref{equation3}) are $x_3x_5$ and $x_4x_6$. Since  $x_3x_6$  divides their  product, by Lemma \ref{Schmitt} it follows that $x_3x_5=x_4x_6=0$.  Thus all generators listed in (\ref{generators2}) vanish at ${\bf x}$.
 \par\smallskip\noindent
\underline{Case 2.2}: Suppose that $x_2\ne0$, so that $x_1x_3-x_2^2=0$. Then the rows of the matrix 
$$\left(\begin{array}{cc}
x_1&x_2\\
x_2&x_3
\end{array}\right)$$
\noindent
are proportional.
\par\smallskip\noindent
\underline{Case 2.2.1}:  Suppose that one of the rows of the above matrix is zero. Then $x_2=0$, so that $x_1x_2x_3=x_2x_4=x_2x_5=x_2x_6=0$.  Hence from (\ref{equation1}) we have that $x_1x_4=0$ and from (\ref{equation4}) we have that $x_1x_5=0$.  From (\ref{equation2}) we then get $x_3x_6=0$. Since $x_3x_6$ divides the product of $x_3x_5$ and $x_4x_6$, which are the summands surviving in (\ref{equation3}), it follows that $x_3x_5=x_4x_6=0$.  Thus all generators listed in (\ref{generators2}) vanish at ${\bf x}$. Now suppose that $x_3=0$, so that $x_3x_5=x_3x_6=0$. Then from (\ref{equation3}) we get $x_4x_6=0$, and from (\ref{equation4}) we get $x_1x_5=0$. Thus all generators listed in (\ref{generators2}) vanish at ${\bf x}$.
\par\smallskip\noindent
\underline{Case 2.2.2}: Suppose that none of the rows of the above matrix is zero. Then there is a nonzero $\lambda\in K$ such that 
\begin{eqnarray}
x_1&=&\lambda x_2=\lambda^2x_3,\nonumber\\
x_2&=&\lambda x_3.\label{lambda}
\end{eqnarray}
\noindent
Replacing these equalities in (\ref{equation2}) and (\ref{equation3}) gives
\begin{eqnarray}
\lambda^2x_3x_4+\lambda x_3x_5+x_3x_6&=&0\label{equation2'}\\
\lambda x_3x_4+x_3x_5+x_4x_6&=&0\label{equation3'}
\end{eqnarray}
\noindent
Comparing (\ref{equation2'}) with $\lambda$(\ref{equation3'}) further yields:
\begin{equation}\label{six}x_3x_6=\lambda x_4x_6.\end{equation}
\smallskip\noindent
\underline{Case 2.2.2.1}: Suppose that $x_6=0$, so that $x_2x_6=x_3x_6=x_4x_6=0$.  Then from (\ref{equation4}) we get $x_1x_5=0$. Since $x_1x_5$ divides the product of $x_1x_4$ and $x_2x_5$, which are the terms surviving in (\ref{equation2}), we conclude that $x_1x_4=x_2x_5=0$. From (\ref{equation1}) we then get that $x_1x_2x_3=0$. Moreover, since $x_2x_5$ divides the product of $x_2x_4$ and $x_3x_5$, which are the terms surviving in (\ref{equation3}), we also conclude that $x_2x_4=x_3x_5=0$. Thus all generators listed in (\ref{generators2}) vanish at ${\bf x}$.
\par\smallskip\noindent
\underline{Case 2.2.2.2}: Suppose that $x_6\ne0$. Then from (\ref{six}) we get
\begin{equation}\label{six2}x_3=\lambda x_4,\end{equation}
\noindent
which, together with (\ref{lambda}), implies that 
\begin{eqnarray}
x_1&=&\lambda^3x_4,\nonumber\\
x_2&=&\lambda^2x_4.\label{26}
\end{eqnarray}
\noindent
Replacing the equalities (\ref{26}) in (\ref{equation4}) gives
\begin{equation}\label{27}\lambda^3x_4x_5+\lambda^2x_4x_6=0,\qquad\mbox{i.e.,}\qquad\lambda x_4x_5+x_4x_6=0.\end{equation}
\noindent
On the other hand, replacing (\ref{six2}) and (\ref{26}) in (\ref{equation2}) we have:
$$\lambda^3x_4^2+\lambda^2x_4x_5+\lambda x_4x_6=0,$$
\noindent
so that, by (\ref{27}), 
$$\lambda^3 x_4^2=0,\qquad\mbox{i.e.,}\qquad x_4=0.$$
\noindent
Thus, by (\ref{lambda}) and (\ref{six2}), we have that $x_1=x_2=x_3=0$. Thus all generators listed in (\ref{generators2}) vanish at ${\bf x}$.
\par\smallskip\noindent
We have examined all possible cases, hence our claim is proven.
}
\end{example}
\begin{remark}{\rm
Four elements generating the ideal $I$ of Example \ref{example4} up to radical had already been found, by a different approach, in \cite{B8}, Example 4. These elements, however, were  strictly depending on the characteristic of $K$. 
}
\end{remark}
\vskip.5truecm
\begin{center}
\textsc{ACKNOWLEDGEMENT}
\end{center}
\noindent
The author is indebted to Marcel Morales for posing the question which motivated this paper.

\end{document}